\documentclass[a4paper,12pt]{article}

\usepackage{amssymb}

\usepackage{amsmath}
\usepackage{amsfonts}
\usepackage{amscd}
\usepackage[english]{babel}
\title{Left-invariant almost nearly K$\mathrm{\ddot{a}}$hler structures on $SU(2)\times SU(2)$ in the tetrahedron visualization
for $\mathbb{CP}^3$.}
\author{Natalia Daurtseva\\ Russia, Kemerovo State University\\
natali0112@ngs.ru}
\date{ }
\begin{document}
\maketitle
\begin{abstract}
The set of maximal non-integrable structures $(SU(2)\times SU(2),B,I)$, where $B$ is Killing-Cartan metric
is described as subset of $\mathbb{CP}^3$. The visualization of complex projective space $\mathbb{CP}^3$ as
tetrahedron which edges and faces are $\mathbb{CP}^1$ and $\mathbb{CP}^2$ is used.
\end{abstract}
{\bf 1. Preliminaries.}

{\bf Definition 1. }{\it An almost complex structure $J$ on  a smooth manifold $M$ is an endomorphism
$$
J:TM\longrightarrow TM
$$
of the tangent bundle such that $J^2=-Id$.}

The almost complex manifold necessarily has even dimension. It is known \cite{1},
that every almost complex structure defines the orientation on $M$. It follows from the fact that for all
$J_x$-independent vectors $X_1,\dots,X_n\in T_xM$ the vectors $X_1,\dots,X_n,J_xX_1,\dots,J_xX_n\in T_xM$
are linear independent and give the same orientation of $T_xM$.

{\bf Definition 2. } {\it Manifold $(M^{2n},g,J_0)$ with almost complex structure $J_0$ and Riemannian metric $g$
is called almost Hermitian, if  $J_0$ is orthogonal with respect to metric $g$, i.e.
$$
g(J_0 X,J_0 Y)=g(X,Y),
$$
for all vector fields $X$ and $Y$ on $M$.}

Let $M$ is 6-dimensional group Lie $G$. In this article we take interest in the left-invariant
orthogonal almost complex structures, giving the same orientation on the group $G=SU(2)\times SU(2)$.
This reduces study of almost complex structure
to reseach its restriction on the Lie algebra $\mathfrak{g}=T_eG$ of $G$. Study of the space of left-invariant
orthogonal almost complex structures on $G$ reduces to reseach  the set $\mathcal{Z}$ of all endomorphisms
$I:\mathfrak{g}\longrightarrow\mathfrak{g}$, preserving orientation, such that $I^2=-1$, and $g_e(IX,IY)=g_e(X,Y)$,
for all $X,Y\in\mathfrak{g}$. Group $SO(6)$ acts transitively on the set $\mathcal{Z}$, with isotropy
$U(3)$, so $\mathcal{Z}$ is homogeneous space $SO(6)/U(3)=\mathbb{CP}^3$.

There is isomorphism between space $\mathcal{Z}$ and $\mathbb{CP}^3$ \cite{2}.
Let $e^1,\dots,e^6$ is dual basis of $\mathfrak{g}^*$, orthonormal with respect to metric on $\mathfrak{g}^*$.
$V=\mathbb{C}^4$ denote the standard representation of $SU(4)$ and let $(v^0,v^1,v^2,v^3)$
be a unitary basis of $V$. Then $\Lambda^2V$ is the complexification of a real vector space which we
identify with $\mathfrak{g}^*$. This identification can be chosen in such a way that:
$$
\begin{array}{ll}
2v^0\wedge v^1=e^1+ie^2 & 2v^2\wedge v^3=e^1-ie^2\\
2v^0\wedge v^2=e^3+ie^4 & 2v^3\wedge v^1=e^3-ie^4\\
2v^0\wedge v^3=e^5+ie^6 & 2v^1\wedge v^2=e^5-ie^6
\end{array}
$$

A point $I\in\mathcal{Z}$ corresponds to a totally isotropic subspace of the complexification of $\mathfrak{g}^*$.
Namely to the $i$-eigenspace of $I$, which consists of vectors $v-iJv$. Any such subspace equals \cite{2}:
$$
V_u=\{u\wedge v:v\in V\}\subset\Lambda^2V,
$$
which depends on $[u]\in P(V)=\mathbb{CP}^3$.
$$
J\longrightarrow (\mathfrak{g}^*)^{1,0}\longrightarrow V_u\longrightarrow [u]\in\mathbb{CP}^3
$$
For example, if $Ie^1=-e^2, Ie^3=-e^4, Ie^5=-e^6$, then $i$-eigenspace of $I$ consists of linear combination of vectors
$$
2v^0\wedge v^1=e^1+ie^2, 2v^0\wedge v^2=e^3+ie^4, 2v^0\wedge v^3=e^5+ie^6
$$

In \cite{2} was offered to visualize $\mathbb{CP}^3$ as a tetrahedron, in which the edges end faces represent projective
subspaces $\mathbb{CP}^1$ and $\mathbb{CP}^2$. The almost complex structures $I_0$, $I_1$, $I_2$, $I_3$:
$$
\begin{array}{c}
I_0e^1=-e^2, I_0e^3=-e^4, I_0e^5=-e^6;\\
I_1e^1=-e^2, I_1e^3=e^4, I_1e^5=e^6;\\
I_2e^1=e^2, I_2e^3=-e^4, I_2e^5=e^6;\\
I_3e^1=e^2, I_3e^3=e^4, I_3e^5=-e^6;\\
\end{array}
$$
are points $[1,0,0,0], [0,1,0,0], [0,0,1,0], [0,0,0,1]\in\mathbb{CP}^3$. Each of these structures associated
with the fundamental 2-forms, by formula:
$$
\begin{array}{c}
\omega(X,Y)=g(IX,Y)\\
\omega_0=e^1\wedge e^2+e^3\wedge e^4+e^5\wedge e^6;\\
\omega_1=e^1\wedge e^2-e^3\wedge e^4-e^5\wedge e^6;\\
\omega_2=-e^1\wedge e^2+e^3\wedge e^4-e^5\wedge e^6;\\
\omega_3=-e^1\wedge e^2-e^3\wedge e^4+e^5\wedge e^6.
\end{array}
$$

Since the metric $g$ is fixed, one may equally well label points of  $\mathcal{Z}$ by their corresponding
fundamental forms.
Coordinates $[v^0,v^1,v^2,v^3]$ on $\mathbb{CP}^3$ are helpful to visualize $\mathcal{Z}$ as solid tetrahedron with
vertices $\omega_0, \omega_1, \omega_2, \omega_3$.

Let $z=[z^0,z^1,z^2,z^3]$ and $u=[u^0,u^1,u^2,u^3]$ are arbitrary
points in $\mathbb{CP}^3$. Any point of the "edge" $E_{zu}$, containing $z$ and $u$ in $\mathbb{CP}^3$
is linear combination $\alpha z+\beta u$, where $\alpha,\beta\in\mathbb{C},\quad |\alpha|+|\beta|\neq 0$.
But $\alpha z+\beta u$ and $c(\alpha z+\beta u)$ define the same point in $\mathbb{CP}^3$, for any $c\in\mathbb{C}^*$.
So "edge",
$E_{zu}=\{[\alpha,\beta]:\alpha,\beta\in\mathbb{C},\quad |\alpha|+|\beta|\neq 0\}=\mathbb{CP}^1=S^2$.
The equatorial circle in the "edge" $E_{zu}$ we denote as $C_{zu}$.

{\bf Lemma 1. }{\it The fundamental 2-form $\omega$ of "edge" $E_{01}$, containing vertices $\omega_0$ and $\omega_1$
has the form:
$$
\omega=e^1\wedge e^2+r(e^3\wedge e^4+e^5\wedge e^6)+u(e^3\wedge e^5+e^6\wedge e^4)+x(e^3\wedge e^6+e^4\wedge e^5)
$$
where $r^2+u^2+x^2=1$}

{\bf Proof:}  For any  $\omega\in E_{01}$ there exist numbers $s,c_1,c_2\in\mathbb{R}$, $s^2+c_1^2+c_2^2=1$
such as corresponding almost complex structure
$I=sI_0+(c_1+ic_2)I_1$.
$$
s[1,0,0,0]+(c_1+ic_2)[0,1,0,0]=[s, c_1+ic_2,0 ,0]
$$
Space $V_{[s, c_1+ic_2,0 ,0]}=\{sv^0+(c_1+ic_2)v^1\wedge u, u\in V\}$.
$$
(sv^0+(c_1+ic_2)v^1)\wedge v^0=(c_1+ic_2)v^1\wedge v^0=-\frac12(c_1+ic_2)(e^1+ie^2)=
$$
$$
=\frac12(-c_1e^1+c_2e^2+i(-c_2e^1-c_1e^2))
$$
i.e.:
$$
I(-c_1e^1+c_2e^2)=c_2e^1+c_1e^2
$$
$$
(sv^0+(c_1+ic_2)v^1)\wedge v^1=sv^0\wedge v^1=\frac12 s(e^1+ie^2)=\frac12(se^1+ise^2)
$$
$$
I(se^1)=-se^2
$$
$$
(sv^0+(c_1+ic_2)v^1)\wedge v^2=sv^0\wedge v^2+(c_1+ic_2)v^1\wedge v^2=\frac12(s(e^3+ie^4)+(c_1+ic_2)(e^5-ie^6))=
$$
$$
\frac12(se^3+c_1e^5+c_2e^6+i(se^4+c_2e^5-c_1e^6))
$$
$$
I(se^3+c_1e^5+c_2e^6)=-(se^4+c_2e^5-c_1e^6)
$$
$$
(sv^0+(c_1+ic_2)v^1)\wedge v^3=sv^0\wedge v^3+(c_1+ic_2)v^1\wedge v^3=\frac12(s(e^5+ie^6)+(c_1+ic_2)(-e^3+ie^4))=
$$
$$
\frac12(se^5-c_1e^3-c_2e^4+i(se^6-c_2e^3+c_1e^4))
$$
$$
I(se^5-c_1e^3-c_2e^4)=-(se^6-c_2e^3+c_1e^4)
$$
Then $\omega\in E_{01}$, corresponding to $I$ has the form:
$$
\omega=(-c_1e^1+c^2e^2)\wedge (-c_2e^1-c^1e^2)+se^1\wedge se^2+(se^3+c_1e^5+c_2e^6)\wedge (se^4+c_2e^5-c_1e^6)+
$$
$$
+(se^5-c_1e^3-c_2e^4)\wedge (se^6-c_2e^3+c_1e^4)=e^1\wedge e^2+(s^2-c_1^2-c_2^2)e^3\wedge e^4+2sc_2e^3\wedge e^5-
$$
$$
-2sc_1e^3\wedge e^6-2sc_1e^4\wedge e^5-2sc_2e^4\wedge e^6+(s^2-c_1^2-c_2^2)e^5\wedge e^6
$$
Let $r=s^2-c_1^2-c_2^2=2s^2-1$, $u=2sc_2$, $x=-2sc_1$, then $r^2+x^2+u^2=1$, and
$$
\omega=\left(
\begin{array}{cccccc}
0 & 1 & 0 & 0 & 0 & 0\\
-1 & 0 & 0 & 0 & 0 & 0\\
0 & 0 & 0 & t & u & -x\\
0 & 0 & -t & 0 & -x & -u\\
0 & 0 & -u & x & 0 & t\\
0 & 0 & x & u & -t & 0\\
\end{array}
\right)
$$
\begin{flushright}$\Box$\end{flushright}

The following notions was defined in the \cite{2}.

{\bf Definition 3. }{\it Let $\sigma=e\wedge f$ - decomposable unit 2-form, the generalized edge $\lceil\sigma\rfloor$
is the set
$$
\lceil\sigma\rfloor=\{\sigma+\tau\in\mathcal{Z}:\tau\in\Lambda^2\langle e,f \rangle^\perp\}.
$$}

{\bf Definition 4. }{\it The polar set of an arbitrary non-zero 2-form $\sigma$ is
$$
\langle\sigma\rangle^{\perp}_{\mathcal{Z}}=\{\omega\in\mathcal{Z}:g(\omega,\sigma)=0\}
$$}

{\bf 2. Group $SU(2)\times SU(2)$.} As Lie group $G$ we will take group $SU(2)\times SU(2)$. Lie algebra
$$
\mathfrak{su}(2)=\{\left(
\begin{array}{cc}
ix_1 & x_2+ix_3\\
x_2-ix_3 & -ix_1
\end{array}
\right):(x_1,x_2,x_3)\in\mathbb{R}^3\}=\mathbb{R}^3,
$$
so $\mathfrak{g}=\mathfrak{su}(2)\times\mathfrak{su}(2)=\mathbb{R}_1^3\times\mathbb{R}_2^3$.
If $e_1, e_2, e_3$ - standard basis $\mathbb{R}_1^3$ and $e_4, e_5, e_6$ - standard basis $\mathbb{R}_2^3$,
then $[e_1,e_2]=e_3$, $[e_1,e_3]=-e_2$, $[e_2,e_3]=e_1$, $[e_4,e_5]=e_6$, $[e_4,e_6]=-e_5$, $[e_5,e_6]=e_4$, $[e_i,e_j]=0$
for $i=1,2,3; j=4,5,6$.

Killing-Cartan form $B(X,Y)=\mbox{tr}(\mbox{ad}X\circ\mbox{ad}Y)$ gives metric
$B(X,Y)=2(X_1,Y_1)+2(X_2,Y_2)$, where $(\ ,\ )$ - standard metric $\mathbb{R}^3$,
$X=(X_1,Y_1),\ Y=(X_2,Y_2)\in\mathbb{R}_1^3\times\mathbb{R}_2^3$
on $SU(2)\times SU(2)$.

The construction of Hopf bundle $\pi:S^3\stackrel{S^1}{\longrightarrow}\mathbb{CP}^1$ may be used for $S^3=SU(2)$.
Therefore, one has for $SU(2)\times SU(2)$:
$$
\pi:SU(2)\times SU(2)\stackrel{SU(1)\times SU(1)}{\longrightarrow}SU(2)/SU(1)\times SU(2)/SU(1)
$$
One can assume that $e_1$ and $e_4$ are tangent at the identity to the first and second factors of fiber $SU(1)\times SU(1)$,
and $\pi_*(e_2)$, $\pi_*(e_3)$, ($\pi_*(e_5)$, $\pi_*(e_6))$ are tangent to first (second) factor of the base.
Construction of the Hopf bundle allows define the canonical 2-form $\omega=e^1\wedge e^4+e^2\wedge e^3+e^5\wedge e^6$
on $SU(2)\times SU(2)$. It is known \cite{1}, that for almost complex structure $I_0$, such as $I_0e_1=e_4$, $I_0e_2=e_3$, $I_0e_5=e_6$,
structure $(B,I_0,\omega)$ is Hermitian, i.e. almost complex structure $I_0$ is integrable.

{\bf Statement 1.}\cite{3}{\it Orthogonal left-invariant almost complex structure $I$ is integrable, i.e.
$(SU(2)\times SU(2),B,I)$ is Hermitian, if and only if corresponding endomorphism is following:
$$
I=\left(\begin{array}{cc}
O_1 & 0\\
0 & O_2
\end{array}\right)I_0
\left(\begin{array}{cc}
O_1^T & 0\\
0 & O_2^T
\end{array}\right),
$$
where $O_1,O_2\in SO(3)$}

Vanishing of Nijenhuis's tensor is criterion of integrability of almost complex structure.
In case of homogeneous space, in particularly
Lie group $SU(2)\times SU(2)$ one can take only restriction of this
tensor on the Lie algebra.
$$
N(X,Y)=[IX,IY]-[X,Y]-I[X,IY]-I[IX,Y],
$$
$X,Y\in\mathfrak{su(2)}\times\mathfrak{su(2)}$. So the set of complex structures is reperesented as
set on which functional of tensor Nijenhuis norm
$$\|N\|^2=\sum_{i,j,k=1}^{6}(N_{ij}^k)^2:I\mapsto \|N_I\|^2$$
is minimal (namely, is equal to zero). What kind of set, given maximum of this functional?

{\bf Statement 2.}
{\it Functional on $\mathcal{Z}$ of tensor Nijenhuis norm on $(SU(2)\times SU(2), B,
I)$, takes maximal values on the following almost complex structures:
$$
I<e_1,e_2,e_3>\subset<e_4,e_5,e_6>;
I<e_4,e_5,e_6>\subset<e_1,e_2,e_3>
$$
}

{\bf Proof: }
Let $I\in\mathcal{Z}$ is invariant almost complex structure.
$I=\left(
\begin{array}{cc}
A & B\\
-B^T & C
\end{array}\right),$
$-A^T=A=\left(
\begin{array}{ccc}
0 & a_1 & a_2\\
-a_1 & 0 & a_3\\
-a_2 & -a_3 & 0
\end{array}
\right)$,
$-C^T=C=\left(
\begin{array}{ccc}
0 & c_1 & c_2\\
-c_1 & 0 & c_3\\
-c_2 & -c_3 & 0
\end{array}
\right)$, $B=\left(
\begin{array}{ccc}
b_1 & b_2 & b_3\\
b_4 & b_5 & b_6\\
b_7 & b_8 & b_9
\end{array}
\right)$.

Tensor Nijenhuis norm is
$$
\frac{1}{24}\|N\|^2=(b_7c_1-c_3b_9+b_6b_1-b_3b_4)^2+(b_5b_9-b_6b_8-b_4a_1-a_2b_7)^2+
$$
$$
+(-b_4b_2+b_1b_5+b_3a_2+b_6a_3)^2+(-a_2b_8-b_4b_9+b_6b_7-b_5a_1)^2+
$$
$$
+(b_2b_9-b_8b_3-b_1a_1+a_3b_7)^2+(-a_3^2-a_2^2-a_1^2+1)^2+(-a_2b_9-b_5b_7+b_4b_8-b_6a_1)^2+
$$
$$
+(b_8c_3+c_2b_7-b_1b_5+b_4b_2)^2+(a_3b_9-b_7b_2+b_1b_8-b_3a_1)^2+
$$
$$
+(a_3b_8-b_9b_1+b_3b_7-b_2a_1)^2+(-c_2b_4-b_5c_3+b_7b_2-b_1b_8)^2+
$$
$$
+(-b_1a_2-b_2b_6+b_5b_3-b_4a_3)^2+(-c_2b_1-b_2c_3+b_4b_8-b_5b_7)^2+
$$
$$
+(-c_3b_6-b_9b_1+b_3b_7+b_4c_1)^2+(c_2b_9+b_2b_6-b_5b_3+c_1b_8)^2+
$$
$$
+(b_1c_1-c_3b_3-b_6b_7+b_4b_9)^2+(c_2b_6-b_2b_9+b_8b_3+c_1b_5)^2+
$$
$$
+(-b_2a_2-b_3b_4+b_6b_1-b_5a_3)^2+(-c_3^2-c^2_2-c^2_1+1)^2+(c_1b_2+c_2b_3-b_6b_8+b_5b_9)^2
$$

Some necessary conditions for  $I^2=-1$ are:
$$
\begin{array}{ll}
 & a_1b_4+a_2b_7-c_1b_2-c_2b_3=0,\\
a_1^2+a_2^2+b_1^2+b^2_2+b_3^2=1, & a_1b_5+a_2b_8+b_1c_1-b_3c_3=0,\\
a_1^2+a_3^2+b_4^2+b_5^2+b_6^2=1, &  a_1b_6+a_2b_9+b_1c_2+b_2c_3=0,\\
a_2^2+a_3^2+b_7^2+b_8^2+b_9^2=1, & -a_1b_1+a_3b_7-c_1b_5-c_2b_6=0,\\
c_1^2+c_2^2+b_1^2+b^2_2+b_3^2=1, & -a_1b_3+a_3b_9+b_4c_2+b_5c_3=0,\\
c_1^2+c_3^2+b_4^2+b_5^2+b_6^2=1, & -a_2b_1-a_3b_4-b_8c_1-b_9c_2=0,\\
c_2^2+c_3^2+b_7^2+b_8^2+b_9^2=1, & -a_2b_2-a_3b_5+b_7c_1-b_9c_3=0,\\
 & -a_2b_3-a_3b_6+b_7c_2+b_8c_3=0\\
 & -a_1b_2+a_3b_8+b_4c_1-b_6c_3=0.
\end{array}
$$
The equation
$$a_1^2+a_2^2+a_3^2=c_1^2+c_2^2+c_3^2$$
follows from the above left six conditions. Let transform
$\frac1{24}\|N\|^2$ taking into consideration above conditions:
$$
\frac{1}{24}\|N\|^2=(b_7c_1-c_3b_9)^2+2(b_7c_1-c_3b_9)(b_6b_1-b_3b_4)+(b_6b_1-b_3b_4)^2+
$$
$$
+(b_4a_1+a_2b_7)^2-2(b_4a_1+a_2b_7)(b_5b_9-b_6b_8)+(b_5b_9-b_6b_8)^2+
$$
$$
+(b_3a_2+b_6a_3)^2-2(b_3a_2+b_6a_3)(b_4b_2-b_1b_5)+(b_4b_2-b_1b_5)^2+
$$
$$
+(b_5a_1+b_8a_2)^2+2(b_5a_1+b_8a_2)(b_4b_9-b_6b_7)+(b_4b_9-b_6b_7)^2+
$$
$$
+(b_7a_3-b_1a_1)^2-2(b_7a_3-b_1a_1)(b_8b_3-b_2b_9)+(b_8b_3-b_2b_9)^2+
$$
$$
2(1-c_1^2-c_2^2-c_3^2)^2+(b_9a_2+b_6a_1)^2-2(b_9a_2+b_6a_1)(b_8b_4-b_5b_7)+(b_8b_4-b_5b_7)^2+
$$
$$
+(b_8c_3+b_7c_2)^2+2(b_8c_3+b_7c_2)(b_4b_2-b_1b_5)+(b_4b_2-b_1b_5)^2+
$$
$$
+(b_9a_3-b_3a_1)^2-2(b_9a_3-b_3a_1)(b_7b_2-b_1b_8)+(b_7b_2-b_1b_8)^2+
$$
$$
+(b_8a_3-b_2a_1)^2+2(b_8a_3-b_2a_1)(b_7b_3-b_1b_9)+(b_7b_3-b_1b_9)^2+
$$
$$
+(b_4c_2+b_5c_3)^2-2(b_4c_2+b_5c_3)(b_7b_2-b_1b_8)+(b_7b_2-b_1b_8)^2+
$$
$$
+(b_1a_2+b_4a_3)^2+2(b_1a_2+b_4a_3)(b_6b_2-b_3b_5)+(b_6b_2-b_3b_5)^2+
$$
$$
+(b_1c_2+b_2c_3)^2-2(b_1c_2+b_2c_3)(b_8b_4-b_5b_7)+(b_8b_4-b_5b_7)^2+
$$
$$
+(b_4c_1-b_6c_3)^2+2(b_4c_1-b_6c_3)(b_7b_3-b_1b_9)+(b_7b_3-b_1b_9)^2+
$$
$$
+(b_9c_2+b_8c_1)^2+2(b_9c_2+b_8c_1)(b_6b_2-b_3b_5)+(b_6b_2-b_3b_5)^2+
$$
$$
+(b_1c_1-b_3c_3)^2+2(b_1c_1-b_3c_3)(b_4b_9-b_6b_7)+(b_4b_9-b_6b_7)^2+
$$
$$
+(b_6c_2+b_5c_1)^2+2(b_6c_2+b_5c_1)(b_8b_3-b_2b_9)+(b_8b_3-b_2b_9)^2+
$$
$$
+(b_2a_2+a_3b_5)^2-2(b_2a_2+a_3b_5)(b_6b_1-b_3b_4)+(b_6b_1-b_3b_4)^2+
$$
$$
+(b_2c_1+c_2b_3)^2+2(b_2c_1+c_2b_3)(b_5b_9-b_6b_8)+(b_5b_9-b_6b_8)^2=
$$
$$
\left[
\begin{array}{ccc}
b_2a_2+a_3b_5=b_7c_1-b_9c_3; & b_4a_1+a_2b_7=b_2c_1+c_2b_3; & b_3a_2+a_3b_6=b_7c_2+b_8c_3;\\
b_5a_1+b_8a_2=-b_1c_1+b_3c_3; & b_7a_3-b_1a_1=b_6c_2+b_5c_1; & b_9a_2+b_6a_1=-b_1c_2-b_2c_3;\\
b_9a_3-b_3a_1=-b_4c_2-b_5c_3; & b_8a_3-b_2a_1=-b_4c_1+b_6c_3; & b_1a_2+b_4a_3=-b_9c_2-b_8c_1.
\end{array}
\right]
$$
$$
=2((b_6b_1-b_3b_4)^2+(b_5b_9-b_6b_8)^2+(b_1b_5-b_4b^2)^2+(b_6b_7-b_4b_9)^2+
$$
$$
+(b_2b_9-b_3b_8)^2+(b_1b_8-b_2b_7)^2+(b_4b_8-b_5b_7)^2+(b_3b_7-b_1b_9)^2+(b_2b_6-b_3b_5)^2+
$$
$$
+(b_7c_1-c_3b_9)^2+(c_1b_2+c_2b_3)^2+(c_2b_7+c_3b_8)^2+(c_1b_1-b_3c_3)^2+(c_1b_5+c_2b_6)^2+
$$
$$
+(c_2b_4+c_3b_5)^2+(c_2b_1+c_3b_2)^2+(c_1b_4-c_3b_6)^2+(c_2b_9+c_1b_8)^2+(1-c_1^2-c_2^2-c_3^2)^2).
$$
The sum of first nine summands is norm of matrix $B^a$, where
$B^a$ is matrix of algebraical complements to $B$. Assume that
$\det B=0$, from the proof of Statement 1 \cite{3} one can see that it is
sufficient condition for integrability, so corresponding
almost complex structure gives $N=0$. We interest in maximal values of
$\|N\|^2$, thus $\det B\neq 0$. We obtain:
$$
\|B^a\|^2=\mbox{tr} (B^{aT}B^a)=\mbox{tr} (\det B\cdot B^{-1}\cdot\det B\cdot
(B^T)^{-1})=(\det B)^2\mbox{tr}(B^TB)^{-1}
$$
Condition $I^2=-1$ gives $B^TB=1+C^2$, so:
$$
\|B^a\|^2=(\det
B)^2\cdot\mbox{tr}(1+C^2)^{-1}=\det(1+C^2)\cdot\mbox{tr}(1+C^2)^{-1}
$$
where $\lambda_1,\ \lambda_2,\ \lambda_3$ are proper values of
matrix $1+C^2$:
$$
\|B^a\|^2=\lambda_1\lambda_2\lambda_3(\lambda_1^{-1}+\lambda_2^{-1}+\lambda_3^{-1})=
\lambda_1\lambda_2+\lambda_2\lambda_3+\lambda_1\lambda_3
$$
Sum $\lambda_1\lambda_2+\lambda_2\lambda_3+\lambda_1\lambda_3$ is
equal to coefficient multiplying $\lambda$ in $\det(1+C^2-\lambda)$:
$$
\|B^a\|^2=c_1^4+c_2^4+c_3^4-4(c_1^2+c_2^2+c_3^2)+3+2c_1^2c_2^2+2c_1^2c_3^2+2c_2^2c_3^2
$$
Sum of last ten summands $\frac{1}{48}\|N\|^2$ is equal to:
$$
1-c_1^4-c_2^4-c_3^4-2c_1^2c_3^2-2c_1^2c_2^2-2c_2^2c_3^2
$$
Thus:
$$
\frac{1}{48}\|N\|^2=4(1-c_1^2-c_2^2-c_3^2),
$$
$$
\|N\|=8\sqrt{3}\sqrt{1-c_1^2-c_2^2-c_3^2}
$$
Then $\max_{J\in{\mathcal Z}}(\|N\|)=8\sqrt{3}$, this value
is reached on almost complex structures $J$, with $C=A=0$.
\begin{flushright}$\Box$\end{flushright}

{\bf Definition 5.\cite{5}} {\it The structure (M,g,J) is called nearly K$\mathrm{\ddot{a}}$hler
if $\nabla_X(\omega)(X,Y)=0$, for all vector fields X,Y on M. The set of all nearly
K$\mathrm{\ddot{a}}$hler structures $(g, J)$ is denoted ${\mathcal NK}$.}

Almost complex structures from the statement 2 satisfy to condition
$$
\nabla_{e_i}(\omega)(e_i,e_j)=0,\qquad\forall i=1,\dots,6,
$$
but for some vectors $X=u+v$, where $u$ is in first summand $\mathfrak{su}(2)$
and $v$ is in second summand one can find vector $Y$:
$$
\nabla_X(\omega)(X,Y)\neq 0,
$$
for example
$$
\nabla_{e_2+e_4}(\omega)(e_2+e_4,e_1)=-1
$$

So, the maximal non-integrable structures $(B,I)$ from  Statement 2
we will call the almost nearly K$\mathrm{\ddot{a}}$hler (${\mathcal ANK}$).

Let realize integrable almost complex structure $I_0$ and ${\mathcal ANK}$ structure
$I_N=\left(
\begin{array}{cc}
0 & -E\\
E & 0
\end{array}\right)$ on $SU(2)\times SU(2)$ as points in $\mathbb{CP}^3$.

{\bf Lemma 2.} $I_0=[1,0,0,-1]\in\mathbb{CP}^3$, $I_N=[1,1,-1,1]\in\mathbb{CP}^3$.

{\bf Proof:} Find eigenspace, corresponding to $I_0$ in $\Lambda^2V$.
$$
\begin{array}{c}
I_0(e^1+ie^4)=i(e^1+ie^4);\\
e^1+ie^4=v^0\wedge v^1+v^2\wedge v^3 + v^0\wedge v^2 - v^3\wedge v^1=(v^0-v^3)\wedge (v^1+v^2);\\
I_0(e^2+ie^3)=i(e^2+ie^3);\\
e^2+ie^3=-iv^0\wedge v^1+iv^2\wedge v^3 + iv^0\wedge v^2 + iv^3\wedge v^1=(v^0-v^3)\wedge (-iv^1+iv^2);\\
I_0(e^5+ie^6)=i(e^5+ie^6);\\
e^5+ie^6=2v^0\wedge v^3=(v^0-v^3)\wedge 2v^3.
\end{array}
$$
Therefore $I_0$ corresponds to  $V_{v^0-v^3}$ and $I_0=[1,0,0,-1]$.

Similarly for $I_N$.
$$
\begin{array}{c}
I_N(e^1-ie^4)=e^4+ie^1=i(e^1-ie^4);\\
e^1-ie^4=v^0\wedge v^1+v^2\wedge v^3 - v^0\wedge v^2 + v^3\wedge v^1=(v^0+v^3)\wedge (v^1-v^2);\\
I_N(e^2-ie^5)=e^5+ie^2=i(e^2-ie^5);\\
e^2-ie^5=-iv^0\wedge v^1+iv^2\wedge v^3 - iv^0\wedge v^3 - iv^1\wedge v^2=(-v^0+v^2)\wedge (iv^1+iv^3);\\
I_N(e^3-ie^6)=e^6+ie^3=i(e^3-ie^6);\\
e^3-ie^6=v^0\wedge v^2+v^3\wedge v^1+v^1\wedge v^2 - v^0\wedge v^3=(v^0+v^1)\wedge (v^2-v^3).
\end{array}
$$
Therefore $I_N$ corresponds to $V_{v^0+v^3-v^2+v^1}$, $I_N=[1,1,-1,1]$.
\begin{flushright}$\Box$\end{flushright}

Find the subset in tetrahedron, corresponding to $\mathcal{ANK}$.

Under the statement 2, the 2-form corresponding to $I\in\mathcal{ANK}$ is following:
$$
\omega=e^4\wedge f^1+e^5\wedge f^2+e^6\wedge f^3,
$$
where $f^i=Ae^i,\quad i=1,2,3$. As $Ie^6\subset\langle e^1,e^2,e^3\rangle$, then 2-form $e^5\wedge e^6$
has  the zero coefficient, i.e.
$\mathcal{ANK}\subset\langle e^5\wedge e^6\rangle^{\perp}_{\mathcal{Z}}$.
Elements of $\langle e^5\wedge e^6\rangle^{\perp}_{\mathcal{Z}}$ are described in \cite{2}: take an arbitrary point
$p_+$ on edge $E_{03}=\lceil e^5\wedge e^6\rfloor$ and arbitrary point $p_-\in E_{12}=\lceil -e^5\wedge e^6\rfloor$.
Connect these points by edge (build the sphere with poles $p_+$ and $p_-$). Then union of all equatorial circles of
edges $E_{p_+p_-}$, $p_+\in E_{03}$, $p_-\in E_{12}$ is $\langle e^5\wedge e^6\rangle^{\perp}_{\mathcal{Z}}$,
so $$\mathcal{ANK}\subset \bigcup_{\begin{array}{c} p_+\in E_{03}\\ p_-\in E_{12} \end{array}}C_{p_+p_-}.$$
Let define the general view of form $\omega\in C_{p_+p_-}$.

{\bf Lemma 3.}{\it Form $\omega\in C_{p_+p_-}$ is:
$$
2\omega=(r_++r_-)e^1\wedge e^2+(r_+-r_-)e^3\wedge e^4+(x_++x_-)e^1\wedge e^3+(x_+-x_-)e^4\wedge e^2+
$$
$$
(u_++u_-)e^1\wedge e^4+(u_+-u_-)e^2\wedge e^3+
$$
$$
+\frac{(x_+t_1-u_+t_2)(r_-+1)+(u_-t_2-x_-t_1)(r_++1)}{\sqrt{(r_-+1)(r_++1)}}e^1\wedge e^5+
$$
$$
+\frac{-(x_+t_2+u_+t_1)(r_-+1)+(u_-t_1+x_-t_2)(r_++1)}{\sqrt{(r_-+1)(r_++1)}}e^1\wedge e^6+
$$
$$
+\frac{(x_+t_2+u_+t_1)(r_-+1)+(u_-t_1+x_-t_2)(r_++1)}{\sqrt{(r_-+1)(r_++1)}}e^2\wedge e^5+
$$
$$
+\frac{(x_+t_1-u_+t_2)(r_-+1)-(u_-t_2-x_-t_1)(r_++1)}{\sqrt{(r_-+1)(r_++1)}}e^2\wedge e^6+
$$
$$
+\frac{-t_2(r_++1)(r_-+1)+u_-(x_+t_1-u_+t_2)+x_-(u_+t_1+x_+t_2)}{\sqrt{(r_-+1)(r_++1)}}e^3\wedge e^5+
$$
$$
+\frac{-t_1(r_++1)(r_-+1)-u_-(x_+t_2+u_+t_1)+x_-(x_+t_1-u_+t_2)}{\sqrt{(r_-+1)(r_++1)}}e^3\wedge e^6+
$$
$$
+\frac{-t_1(r_++1)(r_-+1)+u_-(x_+t_2+u_+t_1)-x_-(x_+t_1-u_+t_2)}{\sqrt{(r_-+1)(r_++1)}}e^4\wedge e^5+
$$
$$
+\frac{t_2(r_++1)(r_-+1)+u_-(x_+t_1-u_+t_2)+x_-(u_+t_1+x_+t_2)}{\sqrt{(r_-+1)(r_++1)}}e^4\wedge e^6
$$
for $r_{\pm}\neq -1$;\\
$2\omega=-e^1\wedge e^2+t_2(e^3\wedge e^5+e^4\wedge e^6)+t_1(e^3\wedge e^6-e^4\wedge e^5)$, for $r_{\pm}=-1$;\\
$2\omega=\frac{r-1}{2}e^1\wedge e^2+t_1\sqrt{\frac{r+1}{2}}(e^1\wedge e^6+e^2\wedge e^5)+\frac{u}{2}(e^1\wedge e^4+e^2\wedge e^3)+
\frac{x}{2}(e^1\wedge e^3-e^2\wedge e^4)-t_2\sqrt{\frac{r+1}{2}}(e^1\wedge e^5-e^2\wedge e^6)+\frac{r+1}{2}e^3\wedge e^4+
\frac{xt_1-ut_2}{\sqrt{2(r+1)}}(e^4\wedge e^6+e^3\wedge e^5)+\frac{ut_1+xt_2}{\sqrt{2(r+1)}}(e^4\wedge e^5-e^3\wedge e^6)$, for $r=r_+\neq -1$, $r_-=-1$.\\
$2\omega=\frac{r-1}{2}e^1\wedge e^2+\frac{u}{2}(e^1\wedge e^4-e^2\wedge e^3)+\frac{x}{2}(e^1\wedge e^3+e^2\wedge e^4)+
t_1\sqrt{\frac{r+1}{2}}(e^1\wedge e^6-e^2\wedge e^5)+t_2\sqrt{\frac{r+1}{2}}(e^1\wedge e^5+e^2\wedge e^6)+
\frac{ut_1+xt_2}{\sqrt{2(r+1)}}(e^3\wedge e^6-e^4\wedge e^5)+\frac{ut_2-xt_1}{\sqrt{2(r+1)}}(e^3\wedge e^5+e^4\wedge e^6)-\frac{1+r}{2}e^3\wedge e^4$,
for $r=r_-\neq -1,\ r_+=-1$.}\\

{\bf Proof:} Similarly lemma 2 one can show, that $p_+\in E_{03}$ and $p_-\in E_{12}$
are of view:
$$
p_+=e^5\wedge e^6+r_+(e^1\wedge e^2+e^3\wedge e^4)+x_+(e^1\wedge e^4+e^2\wedge e^3)+u_+(e^1\wedge e^4+e^2\wedge e^3),
$$
$$
p_-=-e^5\wedge e^6+r_-(e^1\wedge e^2-e^3\wedge e^4)+x_-(e^1\wedge e^4-e^2\wedge e^3)+u_-(e^1\wedge e^4-e^2\wedge e^3),
$$
where $r_+^2+x_+^2+u_+^2=1$ и $r_-^2+x_-^2+u_-^2=1$.

Points in $\mathbb{CP}^3$, corresponding $p_+$ and $p_-$ are:
$$
p_-=\left[0,\sqrt{\frac{r_-+1}{2}},\frac{u_-+ix_-}{\sqrt{2(r_-+1)}},0\right],\quad
p_+=\left[\sqrt{\frac{r_++1}{2}},0,0,\frac{-u_++ix_+}{\sqrt{2(r_++1)}}\right], r_{\pm}\neq 0
$$
$$
p_-=[0,0,1,0];\ p_+=[0,0,0,1], \mbox{при}\ r_{\pm}=-1.
$$
Let consider sphere  - "edge" $E_{p_+p_-}$ now. A point $p$ of this edge is:
$$
p=\left\{\begin{array}{l}
\sqrt{\frac{1-t}{2}}p_-+\frac{t_1+it_2}{\sqrt{2(1-t)}}p_+,\qquad t^2+t_1^2+t_2^2=1, t\neq 1\\
p_+,\qquad t=1, t_1=t_2=0
\end{array}\right.
$$
For $t\neq 1, r_{\pm}\neq -1$ we have: $
2p=(t_1+it_2)\sqrt{\frac{r_++1}{1-t}}v^0+\sqrt{(1-t)(r_-+1)}v^1+\sqrt{\frac{1-t}{1+r_-}}(u_-+ix_-)v^2+
\frac{(-u_++ix_+)(t_1+it_2)}{\sqrt{(r_++1)(1-t)}}v^3.
$
Calculate $p\wedge v^0$:\\
$
2p\wedge v^0=\sqrt{(1-t)(r_-+1)}v^1\wedge v^0+\sqrt{\frac{1-t}{1+r_-}}(u_-+ix_-)v^2\wedge v^0+
\frac{(-u_++ix_+)(t_1+it_2)}{\sqrt{(r_++1)(1-t)}}v^3\wedge v^0=
\frac12[\sqrt{(1-t)(r_-+1)}(-e^1-ie^2)+\sqrt{\frac{1-t}{1+r_-}}(u_-+ix_-)(-e^3-ie^4)+
\frac{(-u_++ix_+)(t_1+it_2)}{\sqrt{(r_++1)(1-t)}}(-e^5-ie^6)]=
\frac12[-\sqrt{(1-t)(r_-+1)}e^1-\sqrt{\frac{1-t}{1+r_-}}u_-e^3+\sqrt{\frac{1-t}{1+r_-}}x_-e^4+
\frac{u_+t_1+x_+t_2}{\sqrt{(r_++1)(1-t)}}e^5+\frac{x_+t_1-u_+t_2}{\sqrt{(r_++1)(1-t))}}e^6+i(-\sqrt{(1-t)(r_-+1)}e^2
-\sqrt{\frac{1-t}{1+r_-}}x_-e^3--\sqrt{\frac{1-t}{1+r_-}}u_-e^4-\frac{x_+t_1-u_+t_2}{\sqrt{(r_++1)(1-t))}}e^5
+\frac{u_+t_1+x_+t_2}{\sqrt{(r_++1)(1-t)}}e^6)]=u_0+iv_0
$
Similarly we can calculate $2p\wedge v^j=u_j+iv_j$, $j=1,2,3$. Then find
$\omega=u_0\wedge v_0+u_1\wedge v_1+u_2\wedge v_2+u_3\wedge v_3$.  For $t=0$ we find form in
equatorial circle $C_{p_+p_-}$.

By analogous we can find points of equatorial circle, for $r_{\pm}=-1$, or $r_+=-1, r_-\neq -1$, $r_-=-1, r_+\neq -1$.
\begin{flushright}$\Box$\end{flushright}

{\bf Theorem.}{\it The set of maximal non-integrable structures
$$\mathcal{ANK}=\bigcup C_{p_+p_-},$$ where
$
p_+=e^5\wedge e^6+r(e^1\wedge e^2+e^3\wedge e^4)+x(e^1\wedge e^4+e^2\wedge e^3)+u(e^1\wedge e^4+e^2\wedge e^3),
$
$
p_-=-e^5\wedge e^6-r(e^1\wedge e^2-e^3\wedge e^4)-x(e^1\wedge e^4-e^2\wedge e^3)+u(e^1\wedge e^4-e^2\wedge e^3),
$
$r^2+x^2+u^2=1$
}

{\bf Proof.} As 2-form $\omega$ corresponding to almost complex structure $I\in\mathcal{ANK}$
is of view:
$$
\omega=e^4\wedge f^1+e^5\wedge f^2+e^6\wedge f^3,\mbox{ where }f^i=Ae^i,\quad i=1,2,3
$$
then in formula for $\omega\in C_{p_+p_-}$ we must claim vanishing of coefficient near the forms
$e^1\wedge e^2, e^1\wedge e^3, e^2\wedge e^3, e^4\wedge e^5, e^4\wedge e^6$.
This requirement gives $r_+=-r_-=r$, $x_+=-x_-=x$, $u_+=u_-=u$.
\begin{flushright}$\Box$\end{flushright}

\end{document}